\documentclass[12pt]{amsart}
\usepackage{latexsym}
\usepackage{amsfonts}
\usepackage{amsmath}
\usepackage{amssymb}
\usepackage{mathtext}
\usepackage[mathscr]{eucal}
\newtheorem{propo}{Proposition}

\newtheorem{lemma}[propo]{Lemma}

\newtheorem{theor}[propo]{Theorem}
\newtheorem{theorem}[propo]{Theorem}

\newcommand{\Ker}{\operatorname{Ker}}

\newcommand{\Irr}{{\mathrm {Irr}}}
\newcommand{\IBR}{{\mathrm {IBr}}}
\newcommand{\IBRL}{{\mathrm {IBr}}_{\ell}}

\newcommand{\Ind}{{\mathrm {Ind}}}

\newcommand{\SL}{{\mathrm {SL}}}

\newcommand{\ZZ}{{\mathbb Z}}

\newcommand{\FF}{{\mathbb F}}

\newcommand{\ta}{\hspace{0.5mm}^{2}\hspace*{-0.2mm}}
\newcommand{\tb}{\hspace{0.5mm}^{3}\hspace*{-0.2mm}}
\begin{document}

\title[Cross characteristic representations of $\ta F_{4}(q)$]
{On the Restriction of Cross Characteristic Representations of $\ta
  F_{4}(q)$ to Proper Subgroups}

\author{Frank Himstedt}
\address{Technische Universit\"at M\"unchen\\Zentrum Mathematik, SB-S-MA\\
Boltzmannstr. 3, 85748 Garching, Germany}
\email{himstedt@ma.tum.de}

\author{Hung Ngoc Nguyen}
\address{Department of Mathematics\\Michigan State University\\
East Lansing, MI 48824, U.S.A.}
\email{hungnguyen@math.msu.edu}

\author{Pham Huu Tiep}
\address{Department of Mathematics\\University of Arizona\\
Tucson, AZ 85721, U.S.A.}
\email{tiep@math.arizona.edu}

\begin{abstract}
We prove that the restriction of any
nontrivial representation of the Ree groups $\ta
F_{4}(q), q=2^{2n+1}\geq8$ in odd characteristic to any proper
subgroup is reducible.
We also determine all triples $(K, V, H)$ such that
$K \in \{\ta F_4(2), \ta F_4(2)'\}$, $H$ is a proper subgroup of $K$,
and $V$ is a representation of $K$ in odd
characteristic restricting absolutely irreducibly to $H$.
\end{abstract}

\thanks{Part of this work was done while the authors
were participating in the program on Representation Theory of Finite
Groups and Related Topics at the Mathematical Sciences Research
Institute (MSRI), Berkeley. It is a pleasure to thank the organizers
Professors J. L.~Alperin, M.~Brou\'e, J. F.~Carlson, A. S.~Kleshchev,
J.~Rickard, B.~Srinivasan for generous hospitality and support and
stimulating environment.}
\thanks{The third author gratefully acknowledges the
support of the NSF (grants DMS-0600967 and DMS-0901241).}

\subjclass{Primary 20C20; Secondary 20C33, 20C15}
\keywords{Cross characteristic representations, Irreducible restrictions}

\maketitle

\section{Introduction} \label{sec:intro}

Let $G$ be a quasi-simple group. The problem of classifying all
pairs $(\Phi,H)$, where $\Phi$ is a representation of
$G$ which is absolutely irreducible over a proper subgroup $H$ of $G$, has been
studied recently in many papers. The solution of this problem is a
major part of the Aschbacher-Scott program
to classify maximal subgroups of finite
classical groups (\cite{A}, \cite{Sc}).

When $G$ is a cover of a sporadic group, the problem is largely
computational. The case when $G$ is a cover of symmetric or
alternating groups is quite complicated and almost completed in
\cite{BK}, \cite{KS}, \cite{KT1}, and \cite{Sx}. Our main focus is
on the case where $G$ is a finite group of Lie type. This has been
solved recently in \cite{KT2} when $G$ is of type $A$ and in
\cite{N}, \cite{NTH} when $G$ is of type $G_2$, $\ta G_2$, $\ta B_2$, $\tb
D_4$. In this paper, we solve the problem for the Ree groups
$\ta F_4(q)$. The main result is the following:

\begin{theor}\label{thm:main}
{Let $G = \ta F_4(q)$, $q=2^{2n+1}$, $n\ge 1$ and let $\Phi$ be
  any nontrivial representation of $G$ in characteristic
  $\ell \neq 2$. If $H$ is any proper subgroup of $G$,
  then $\Phi|_H$ is reducible.}
\end{theor}

The structure of the paper is the following. In the next section, we
show that the restriction of any irreducible representation of $\ta
F_4(q)$ in odd characteristic to most of its maximal subgroups is
reducible except possibly for the two maximal parabolic subgroups.
This is done by
using the results on the maximal subgroups \cite{M} and the lower bound on the
degrees of nontrivial irreducible representations of $\ta F_4(q)$
proved in \cite{T}. The small groups $\ta F_4(2)$ as well as $\ta
F_4(2)'$ are handled in $\S3$. Finally, in $\S4$, we finish the
proof of Theorem \ref{thm:main} by showing that the restrictions to
maximal parabolic subgroups are reducible.

\section{Basic Reduction} \label{sec:basreduct}

Let $\FF$ be an algebraically closed field of characteristic
$\ell\neq 2$. Given a finite group $X$, we denote by
$\mathfrak{d}_{\ell}(X)$ and $\mathfrak{m}_{\ell}(X)$ the smallest
and largest dimensions of non-trivial irreducible
$\mathbb{F}X$-modules. As usual, $\Irr(X)$ (resp. $\IBRL(X)$) will
be the set of irreducible  complex characters (resp. irreducible
$\ell$-Brauer characters) of $X$. If $\chi$ is a complex character
of $X$, we denote by $\widehat{\chi}$ the restriction of $\chi$ to
$\ell$-regular elements of $X$.

First we record the following obvious observation:

\begin{lemma}\label{lemma1}
{Let $K$ be a finite group. Suppose that $V$ is an irreducible
$\FF K$-module of dimension greater than one, and $H$ is a proper
subgroup of $K$ such that the restriction $V|_H$ is irreducible.
Then
$$\sqrt{|H/Z(H)|} \geq \mathfrak{m}_{\mathbb{C}}(H)\geq
\mathfrak{m}_{\ell}(H)\geq \dim(V)\geq \mathfrak{d}_{\ell}(K).$$}
\end{lemma}

\begin{lemma}\cite[Corollary (11.29)]{I}\label{lemmaIsaacs}
{Let $K$ be a finite group and $H$ be a normal subgroup of $K$.
Let $\chi\in \Irr(K)$ and $\theta\in\Irr(H)$ be a constituent of
$\chi_H$. Then $\chi(1)/\theta(1)$ divides $|K/H|$.}
\end{lemma}

\begin{lemma}\label{lemmaFeit}
{Let $K$ be a finite group and $H\lhd K$. Suppose $|K/H|=p$ is
prime and $\chi\in \IBR_{\ell}(K)$. Then either
\begin{enumerate}
  \item[(i)] $\chi|_H$ is irreducible or
  \item[(ii)] $\chi|_H=\sum_{i=1}^p\theta_i$, where the $\theta_i$'s are distinct and irreducible.
\end{enumerate}}
\end{lemma}

\begin{proof}
Since $|K/H|$ is prime, the inertia subgroup of an irreducible
constituent of $\chi|_H$ is $H$ or $K$. Now, the claim follows from
Clifford theory, see \cite[III.2.12, III.2.5, III.2.11, III.2.14]{F}.
\end{proof}

In the following theorem, we use the result and notation in \cite{M}
except that the Ree groups will be denoted by $\ta F_{4}(q)$ instead
of $\ta F_{4}(q^2)$.

\begin{theorem}[Reduction Theorem]\label{reduction}
{Let $G = \ta F_{4}(q)$, $q=2^{2n+1}, n\geq1$ and let $\Phi$ be
an irreducible representation of $G$ in characteristic $\ell\neq 2$.
Suppose that $\deg(\Phi)>1$ and $M$ is a maximal subgroup of $G$
such that $\Phi|_{M}$ is irreducible. Then $M$ is $G$-conjugate to
one of the following groups:
\begin{enumerate}
\item[(i)] $P_a=[q^{11}]:(L_2(q)\times (q-1))$, the short-root maximal
parabolic subgroup of order $q^{12}(q-1)(q^2-1)$,
\item[(ii)] $P_b=[q^{10}]:(Sz(q)\times (q-1))$, the long-root maximal
parabolic subgroup of order $q^{12}(q-1)^2(q^{2}+1)$.
\end{enumerate}}
\end{theorem}

\begin{proof}
By \cite[Theorem 1.4]{T}, we have
$\mathfrak{d}_{\ell}(G)=\mathfrak{d}_{\ell}(\ta F_4(q))\geq
(q^4+q^3+q)(q-1)\sqrt{q/2}$ for every $\ell\neq 2$. Next, according
to \cite{M}, if $M$ is a maximal subgroup of $G$, but $M$ is not a
maximal parabolic subgroup, then $M$ is $G$-conjugate to one of the
following groups:
\begin{enumerate}
\item[1)] ${N}_G(\langle t_4\rangle)=C_G(t_4):2\simeq SU_3(q):2$,
\item[2)] $N_G(T_8)\simeq (\mathbb{Z}_{q+1}\times\mathbb{Z}_{q+1}):GL_2(3)$,
\item[3)] $N_G(T_6)\simeq (\mathbb{Z}_{q-\sqrt{2q}+1}\times\mathbb{Z}_{q-\sqrt{2q}+1}):[96]$ if $q>8$,
\item[4)] $N_G(T_7)\simeq (\mathbb{Z}_{q+\sqrt{2q}+1}\times\mathbb{Z}_{q+\sqrt{2q}+1}):[96]$,
\item[5)] $N_G(T_{10})\simeq\mathbb{Z}_{q^2-\sqrt{2q^3}+q-\sqrt{2q}+1}:12$,
\item[6)] $N_G(T_{11})\simeq\mathbb{Z}_{q^2+\sqrt{2q^3}+q+\sqrt{2q}+1}:12$,
\item[7)] $PGU_3(q):2$,
\item[8)] $Sz(q)\wr 2$,
\item[9)] $Sp_4(q):2$,
\item[10)] $\ta F_4(q_0)$ with $q=q_0^\alpha$, $\alpha$ prime.
\end{enumerate}

First, if $M$ is one of the groups from 2) to 6), then by applying Lemma
\ref{lemmaIsaacs}, we have $\mathfrak{m}_{\mathbb{C}}(M)\leq
\max\{|GL_2(3)|,96,12\}<(q^4+q^3+q)(q-1)\sqrt{q/2}\leq
\mathfrak{d}_{\ell}(G)$ for every $q\geq 8$, contradicting Lemma
\ref{lemma1}. Second, we consider the case when $M$ is one of the groups
in 1) and 7). Note that $|SU_3(q)|=|PGU_3(q)|=q^3(q^3+1)(q^2-1)$.
Therefore $\mathfrak{m}_{\mathbb{C}}(M)\leq
\sqrt{2q^3(q^3+1)(q^2-1)}$, which is less than
$(q^4+q^3+q)(q-1)\sqrt{q/2}$ when $q\geq 8$. This leads to a
contradiction again by Lemma \ref{lemma1}. The groups $Sz(q)\wr 2$
and $Sp_4(q):2$ can be excluded similarly since
$|Sz(q)|=q^2(q^2+1)(q-1)$ and $|Sp_4(q)|=q^4(q^2-1)(q^4-1)$.
Finally, assume that $M=\ta F_4(q_0)$ with $q=q_0^\alpha$, $\alpha$
prime. Since $\alpha\mid (2n+1)$, $\alpha\geq 3$. Then we have
$|M|=q_0^{12}(q_0-1)(q_0^3+1)(q_0^4-1)(q_0^6+1)<q_0^{26}\leq
q^{26/3}$. It follows that
$\mathfrak{m}_{\mathbb{C}}(M)<q^{13/3}<(q^4+q^3+q)(q-1)\sqrt{q/2}\leq
\mathfrak{d}_{\ell}(G)$, a contradiction. We have shown that the
restriction of $\Phi$ to any of the groups from 1) to 10) is
reducible, as desired.
\end{proof}

\section{Irreducible Restrictions for $\ta F_4(2)'$ and $\ta F_4(2)$}
\label{sec:resnongeneric}

In the following theorem and its proof, we mainly use the notation
in \cite{Atl1} and \cite{Atl2}.
By $Cl_X$ we mean the conjugacy class $Cl$ in a group $X$.

\begin{theorem}
{Let $G\in\{\ta F_4(2)',\ta F_4(2)\}$ and $\varphi\in
\IBR_\ell(G)$, $\ell\neq 2$. Suppose that $\varphi$ is nonlinear and
$M$ is a proper subgroup of $G$. Then we have
\begin{enumerate}
\item[(i)] When $G=\ta F_4(2)'$, $\varphi|_M$ is irreducible if and only
if one of the following holds:
\begin{enumerate}
\item[(a)] $M=L_3(3)$ or $L_3(3):2$, $\varphi$ is any of the two characters of
degree $26$ when $\ell\neq2,3$ or $\varphi$ is any of the two
characters of degree $27$ when $\ell\neq2,13$;
\item[(b)] $M=L_2(25)$, $\varphi$ is any of the two characters of
degree $26$ when $\ell\neq2,3,5$.
\end{enumerate}
\item[(ii)] When $G=\ta F_4(2)$, $\varphi|_M$ is irreducible if and only
if one of the following holds:
\begin{enumerate}
\item[(a)] $M=\ta F_4(2)'$, $\varphi$ is any character which is an
  extension of an irreducible character of $\ta F_4(2)'$;
\item[(b)] $M=L_3(3)$ or $L_3(3):2$, $\varphi$ is any of the four
irreducible characters of degree $27$ when $\ell\neq2,13$.
\end{enumerate}
\end{enumerate}}
\end{theorem}

\begin{proof}
(i) Every proper subgroup is contained in a maximal one. So we now
temporarily suppose that $M$ is maximal. If $(M,\varphi)$ is a pair
such that $\varphi|_M$ is irreducible, then we will go further by
restricting $\varphi$ to proper subgroups of $M$. Actually, we will
see later that all examples happen when $M$ is maximal. By
\cite{Atl1}, if $M$ is a maximal subgroup of $\ta F_4(2)'$, then $M$
is $\ta F_4(2)'$-conjugate to one of the following subgroups:
\begin{enumerate}
\item[1)] $L_3(3):2$, two non-conjugate subgroups,
\item[2)] $2.[2^8].5.4$,
\item[3)] $L_2(25)$,
\item[4)] $2^2.[2^8].S_3$,
\item[5)] $A_6.2^2$, two non-conjugate subgroups,
\item[6)] $5^2:4A_4$.
\end{enumerate}

According to both complex and modular character tables of $\ta
F_4(2)'$ from \cite{Atl1} and \cite{Atl2}, $\ta F_4(2)'$ has two
irreducible complex characters of degree $26$, which are denoted by
$\chi_2$ and $\chi_3=\overline{\chi_2}$; and two irreducible complex
characters of degree $27$, which are denoted by $\chi_4$ and
$\chi_5=\overline{\chi_4}$. Moreover, any reduction modulo
$\ell\neq2$ of these characters is still irreducible. We also note
that $\ta F_4(2)'$ has an irreducible complex character $\chi_6$ of
degree 78 and $\widehat{\chi_6}\in\IBR_\ell(G)$ for $\ell\neq2,3$.
When $\ell=3$, $\widehat{\chi_6}-\widehat{1_G}$ is irreducible.

1) $M=L_3(3):2$. Note that $\mathfrak{m}_\mathbb{C}(M)\leq52$.
Therefore if $\varphi|_M$ is irreducible then $\varphi(1)=26,27$.
First, if $\varphi(1)=26$ then $\varphi=\widehat{\chi_2}$ or
$\widehat{\chi_3}$. Now we will show that
$\chi_2|_{L_3(3)}=\chi_3|_{L_3(3)}=\chi_8$ where $\chi_8$ is a
character of degree 26 of $L_3(3)$ as denoted in \cite[p. 13]{Atl1}.
Since the value of any irreducible character (different from
$\chi_8$) of $L_3(3)$ of degree $\leq26$ at the class $6A_{L_3(3)}$
is non-negative while $\chi_2(6A_G)=\chi_3(6A_G)=-1$, it follows
that $\chi_2|_{L_3(3)}=\chi_3|_{L_3(3)}=\chi_8$. Furthermore, since
$\widehat{\chi_8}$ is irreducible for every $\ell\neq2,3$,
$\widehat{\chi_2}|_{L_3(3)}$ as well as $\widehat{\chi_3}|_{L_3(3)}$
are irreducible when $\ell\neq2,3$. When $\ell=3$,
$\widehat{\chi_2}|_{M}$ is reducible since $\IBR_3(M)$ does not
contain any irreducible character of degree 26.

Next, we consider $\varphi(1)=27$. Then $\varphi=\widehat{\chi_4}$
or $\widehat{\chi_5}$. Now we will show that
$\chi_4|_{L_3(3)}=\chi_5|_{L_3(3)}=\chi_{11}$, where $\chi_{11}$ is
the unique irreducible character of degree 27 of $L_3(3)$. Note that
$L_3(3)$ has a unique conjugacy class of elements of order 4, which
is denoted by $4A_{L_3(3)}$ and $G$ has 3 classes $4A_G, 4B_G$, and
$4C_G$. Since the value of any character of degree $\leq26$ of
$L_3(3)$ at the class $4A_{L_3(3)}$ is non-negative while
$\chi_2(4A_G)=\chi_2(4B_G)=-2$, it follows that
$4A_{L_3(3)}\subset4C_G$. Now we assume that
$\chi_4|_{L_3(3)}\neq\chi_{11}$. Then $\chi_4|_{L_3(3)}$ is sum of
characters of degrees less than 27 and therefore
$\chi_4|_{L_3(3)}(4A_{L_3(3)})\geq0$, which leads to a contradiction
since $\chi_4(4C_G)=-1$. So
$\chi_4|_{L_3(3)}=\chi_5|_{L_3(3)}=\chi_{11}$ since $\chi_{11}$ is
real. By inspecting the Brauer character table of $L_3(3)$, we see
that $\widehat{\chi_{11}}\in\IBR_\ell(L_3(3))$ for every
$\ell\neq2,13$. Note that $L_3(3)$ does not have any irreducible
13-Brauer character of degree 27. Therefore, when $\ell=13$,
$\widehat{\chi_4}|_{L_3(3)}$ is reducible and so is
$\widehat{\chi_4}|_{M}$.

Now suppose that $M$ is a proper subgroup of $L_3(3)$ or $L_3(3):2$
but $M\neq L_3(3)$, and $\varphi=\widehat{\chi_i}$ for $i=2,3,4,5$.
By \cite[p. 13]{Atl1}, $|M|\leq 432$. It follows that
$\mathfrak{m}_\mathbb{C}(M)\leq\sqrt{|M|}<21$. This implies that
$\varphi|_M$ is reducible since $\widehat{\chi_i}(1)=26$ or $27$.

2) $M=2.[2^8].5.4$. Since
$\mathfrak{m}_\mathbb{C}(M)\leq\sqrt{|M|}<102$, if $\varphi|_M$ is
irreducible then $\varphi(1)<102$ and therefore $\varphi$ must be
the reduction modulo $\ell$ of one of the characters of degrees 26,
27, 78 or $\varphi=\widehat{\chi_6}-\widehat{1_G}$ when $\ell=3$. Note that
$M$ is solvable. Hence, by the Fong-Swan Theorem $\varphi$ lifts to an
irreducible complex representation of $M$, and so
$\varphi(1) \in \{26, 27, 78, 77\}$ divides $|M|$, a contradiction.

3) $M=L_2(25)$. According to \cite[p. 17]{Atl1},
$\mathfrak{m}_\mathbb{C}(M)=26$. Therefore if $\varphi|_M$ is
irreducible then $\varphi=\widehat{\chi_2}$ or $\widehat{\chi_3}$.
First we show that $\chi_2|_M$ and $\chi_3|_M$ are indeed
irreducible. Note that $M$ has a unique class of elements of order
3, which is denoted by $3A_M$, and this class is contained in the
class $3A_G$ of $G$. The value of any character of $M$ of degree
$<26$ at the class $3A_M$ is non-negative. On the other hand
$\chi_2(3A_G)=\chi_3(3A_G)=-1$. It follows that $\chi_2|_M$ and
$\chi_3|_M$ are irreducible. In fact, it is easy to see that
$\chi_2|_M=\chi_3|_M=\chi_{13}$, where $\chi_{13}$ is an irreducible
character of $M$ mentioned in \cite[p. 17]{Atl1}. By inspecting the
Brauer character table of $L_2(25)$, we have that
$\widehat{\chi_{13}}$ is irreducible if and only if $\ell\neq3,5$.

Now suppose that $M$ is a proper subgroup of $L_2(25)$ and
$\varphi=\widehat{\chi_i}$ for $i=2,3$. Then, by \cite[p. 16]{Atl1},
$\mathfrak{m}_\mathbb{C}(M)\leq\sqrt{|M|}<18$. This implies that
$\varphi|_M$ is reducible since $\chi_i(1)=26$.

4) $M=2^2.[2^8].S_3$. This case is treated similarly as the case 2).
We also have that $\varphi|_M$ is reducible for every
$\varphi\in\IBR_\ell(G)$ with $\ell\neq2$.


5) $M=A_6.2^2$ or 6) $M=5^2:4A_4$. We have
$\mathfrak{m}_\mathbb{C}(M)\leq\sqrt{|M|}<38$. Therefore if
$\varphi|_M$ is irreducible then $\varphi(1)=26$ or 27. That means
$\varphi=\widehat{\chi_i}$ for $i=2,3,4,5$. Since $26,27\nmid|M|$,
it follows that $\chi_i|_M$ is reducible and so is
$\widehat{\chi_i}|_M$.

\smallskip

(ii) By \cite{Atl1} and \cite{M}, if $M$ is a maximal subgroup of $\ta F_4(2)$
but $M\neq\ta F_4(2)'$, then $M$ is $\ta F_4(2)$-conjugate to one of
the following subgroups:
\begin{enumerate}
\item[1)] $13:12$,
\item[2)] $2.[2^9].5.4$,
\item[3)] $L_2(25).2$,
\item[4)] $2^2.[2^9].S_3$,
\item[5)] $5^2:4S_4$,
\item[6)] $SU_{3}(2):2$.
\end{enumerate}

1) $M=13:12$. We have
$\mathfrak{d}_\ell(G)\geq27>\mathfrak{m}_\mathbb{C}(M)$ for every
$\ell\neq2$. Therefore $\varphi|_M$ is reducible for any
$\varphi\in\IBR_\ell(G)$.

2) -- 6).
Suppose that $\varphi|_M$ is irreducible. Then $\varphi(1)\leq
\mathfrak{m}_\mathbb{C}(M)\leq\sqrt{|M|}<144$. Moreover, since $M$
is solvable, $\varphi|_M$ is liftable to a complex character of $M$
and therefore $\varphi(1)\mid|M|$. Inspecting both the complex and
modular character tables of $\ta F_4(2)$, we see that $\ta F_4(2)$
does not have any irreducible character satisfying these conditions.

3) $M=L_2(25).2$. We have $\mathfrak{m}_\mathbb{C}(M)\leq52$.
Therefore, if $\varphi|_M$ is irreducible then $\varphi$ is the
reduction modulo $\ell$ of the unique irreducible character of
degree 52 (which we denote by $\chi$) or four irreducible characters
of degree 27 (which we denote by $\chi_4'$, $\chi_4''$, $\chi_5'$
and $\chi_5''$). Here we notice that $\chi_2$ and $\chi_3$ fuse in
$G$ and give the character $\chi$ while $\chi_4$ extends to
$\chi_4'$ and $\chi_4''$ and $\chi_5$ extends to $\chi_5'$ and
$\chi_5''$.

Since $27\nmid|M|$, it follows that $\chi_4'|_M$, $\chi_4''|_M$,
$\chi_5'|_M$ and $\chi_5''|_M$ are reducible. Now it remains to
consider $\varphi=\widehat{\chi}$. As we have mentioned in (i)3),
$\chi_2|_{L_2(25)}=\chi_3|_{L_2(25)}=\chi_{13}$. It follows that
$\chi|_{L_2(25)}=2\chi_{13}$. This and Lemma \ref{lemmaFeit} imply
that $\chi|_M$ is reducible and so is $\widehat{\chi}|_M$.

We have shown that $\varphi|_M$ is reducible if $M$ is not a
subgroup of $\ta F_4(2)'$. Now statement (ii) of the theorem follows
immediately from part (i).
\end{proof}

\section{Restriction to Maximal Parabolic Subgroups}
\label{sec:resparab}


\begin{lemma}\label{la:paraPa}
{Let $P_a$ denote the short-root maximal parabolic subgroup of
order $q^{12}(q^2-1)(q-1)$ of $G = \ta F_{4}(q)$, $q=2^{2n+1}$, $n
\ge 1$ and let $U_a := O_{2}(P_a)$. Then

\begin{enumerate}
\item[(i)] For any odd prime $r$, $O_{r}(P_a) = 1$.

\item[(ii)] The center $Z := Z(U_a)$ of $U_a$ is an
elementary-abelian group of order $q^2$.

\item[(iii)] Let $\lambda \in \IBR_\ell(Z)$, $\ell \neq 2$, be a
non-trivial irreducible Brauer character of $Z$ and~$I$ the inertia
subgroup of $\lambda$ in $P_a$. Then $I$ is solvable.

\item[(iv)] Let $\varphi \in \IBRL(P_a)$, $\ell \neq 2$, be a
faithful irreducible Brauer character of $P_a$. Then~$\varphi$ lifts
to a faithful complex character $\chi$ of $P_a$.
\end{enumerate} }
\end{lemma}

\begin{proof}
(i) Since $O_{r}(P_a), U_a \lhd P_a$ and $O_{r}(P_a) \cap U_a = 1$,
any element $g \in O_{r}(P_a)$ is centralized by $U_a$, which has
order $q^{11}$. Thus $q^{11}$ divides $|C_{P_a}(g)|$. Assuming $g
\neq 1$, we see by \cite[Table A.3]{HH2} that $g$ is $P_a$-conjugate
to a long-root element $\alpha_{12}(1)$.
But then $g$ is a $2$-element, a contradiction. Hence $O_{r}(P_a) = 1$.

(ii) By \cite[(3.6)]{Sh}, we have $Z(U_a) = \{\alpha_{11}(d_1)
\alpha_{12}(d_{12}) \, | \, d_1, d_2 \in \FF_{q} \}$, which is an
elementary-abelian group of order $q^2$.

(iii) From the orders of the centralizers in \cite[Table A.3]{HH2}),
we see that $P_a$ acts transitively on $Z \setminus \{1\}$ and so,
$P_a$ also acts transitively on $\Omega := \IBR_\ell(Z) \setminus
\{1_Z\}$. Hence, $|P_a:I| = q^2-1$ and $U_a \subseteq I$. Consider
$\bar{I} := I / U_a$ which we identify with a subgroup of the Levi
complement $L_a = Z(L_a) \times [L_a, L_a] \cong \ZZ_{q-1} \times
\SL_2(q)$ of $P_a$. The subgroup $\bar{I}_1 := \bar{I} \cap [L_a,
L_a]$ is a normal subgroup of $\bar{I}$ and $\bar{I} / \bar{I}_1$ is
cyclic. So, it is enough to show that $\bar{I}_1$ is solvable.

Since $[L_a, L_a] \lhd L_a$ of index $q-1$ and $L_a$ acts
transitively on $\Omega$, all orbits of $[L_a, L_a]$ on $\Omega$
have the same size, say $t$, and $t$ divides $q^2-1$ and $t$ is a
multiple of $q+1$. Hence, the index of $\bar{I}_1$ in $[L_a, L_a]$
divides $q^2-1$ and is a multiple of $q+1$. So, $\bar{I}_1$ is
isomorphic to a subgroup of $\SL_2(q)$ and $q \, | \, |\bar{I}_1|$
and $|\bar{I}_1| \, | \, q(q-1)$. It is easy to see (for example
using the Bruhat decomposition) that every subgroup of $\SL_2(q)$
containing a Sylow $2$-subgroup and whose order divides $q(q-1)$ is
contained in some Borel subgroup of $\SL_2(q)$. So, $\bar{I}_1$ is
solvable and then also $I$ is solvable.

(iv) Let $\lambda$ be an irreducible constituent of $\varphi|_Z$,
and let $I$ be the inertia subgroup of $\lambda$ in $P_a$. By
Clifford theory, $\varphi = \Ind^{P_a}_I(\psi)$ for some $\psi \in
\IBR_\ell(I)$ whose restriction to $Z$ contains $\lambda$. Since
$\varphi$ is faithful, $\lambda$ is not the trivial character. So
(iii) implies that $I$ is solvable. By the Fong-Swan Theorem, $\psi$
lifts to a complex character $\rho$ of $I$. Hence, $\varphi$ lifts
to the complex character $\chi := \Ind^{P_a}_I(\rho)$.

Assume that $K := \Ker(\chi)$ is non-trivial; in particular $\ell
\neq 0$. If $K$ is not an $\ell$-group, then $K$ contains a
non-trivial $\ell'$-element $g$. Since $\varphi(g) = \chi(g) =
\chi(1) = \varphi(1)$, we see that $\varphi$ is not faithful, a
contradiction. Hence $K$ is an $\ell$-group, and so $O_\ell(P_a)
\neq 1$, contradicting (i).
\end{proof}

\begin{theor}\label{thm:resparab}
{Let $M$ be a maximal parabolic subgroup of $G = \ta F_4(q)$,
  $q=2^{2n+1}$, $n \ge 1$ and let $\varphi \in \IBR_\ell(G)$, $\ell
  \neq 2$, be of degree $>1$. Then $\varphi|_{M}$ is reducible.}
\end{theor}

\begin{proof}
First suppose that $M=P_b$, the long-root parabolic subgroup of $G$.
Then the statement follows from Theorem 1.6 of \cite{T}. So we may
assume that $M=P_a$, the other maximal parabolic subgroup of $G$.
Also assume the contrary: $\varphi|_{P_a}$ is irreducible.

We will consider the two long-root elements $u := \alpha_{12}(1)$ and
$v := \alpha_8(1)$ of $P_a$, in the notation of \cite{HH1}, \cite{HH2}.
Certainly, these two long-root elements are conjugate
in~$G$, so $\varphi(u) = \varphi(v)$. By Lemma \ref{la:paraPa},
$\varphi|_{P_a}$ lifts to a complex irreducible character $\chi$ of
$P_a$ which is faithful. Since $u$ and $v$ are $\ell'$-elements, we
have $\varphi(u) = \chi(u)$ and $\varphi(v) = \chi(v)$. It follows
that
\begin{equation}\label{eq:uv}
  \chi(u) = \chi(v).
\end{equation}
According to \cite{HH2}, the faithful character $\chi$ must be one
of ${_{P_a}\chi}_{36}$, ${_{P_a}\chi}_{37}$ or ${_{P_a}\chi}_{j}(k)$
with $j \in \{34,35,38,39,40\}$. The values of these characters on
$u$ and $v$ are computed in \cite{HH2} and displayed in the following
Table:

\begin{center}
\begin{small}
\begin{tabular}{|c|c|c||c|c|c|}
\hline
\rule{0cm}{0.4cm}
& $u$ & $v$ && $u$ & $v$
\rule[-0.2cm]{0cm}{0.3cm}\\
\hline
\rule{-0.12cm}{0.55cm}
${_{P_a}\chi}_{34}(k)$ & $-\sqrt{\frac{q^7}{2}}$      &
$-(q-1)\sqrt{\frac{q^5}{2}}$ & ${_{P_a}\chi}_{38}(k)$ &
$-(q-1)\sqrt{\frac{q^7}{2}}$ & $-(q-1)^2\sqrt{\frac{q^5}{2}}$
\rule[-0.2cm]{0cm}{0.6cm}\\
\rule{-0.12cm}{0.6cm}
${_{P_a}\chi}_{35}(k)$ & $-\sqrt{\frac{q^7}{2}}$      &
$-(q-1)\sqrt{\frac{q^5}{2}}$ & ${_{P_a}\chi}_{39}(k)$ &
$-(q-1)\sqrt{\frac{q^7}{2}}$ & $-(q-1)^2\sqrt{\frac{q^5}{2}}$
\rule[-0.2cm]{0cm}{0.6cm}\\
\rule{-0.12cm}{0.6cm}
${_{P_a}\chi}_{36}$    & $-(q-1)\sqrt{\frac{q^7}{2}}$ &
$-(q-1)^2\sqrt{\frac{q^5}{2}}$ & ${_{P_a}\chi}_{40}(k)$ & $-q^4(q-1)$
& $q^3(q-1)$
\rule[-0.2cm]{0cm}{0.6cm}\\
\rule{-0.12cm}{0.6cm}
${_{P_a}\chi}_{37}$    & $-(q-1)\sqrt{\frac{q^7}{2}}$ &
$-(q-1)^2\sqrt{\frac{q^5}{2}}$ &&&
\rule[0cm]{0cm}{0.6cm}\\
\hline
\end{tabular}
\end{small}
\end{center}

Now one sees that (\ref{eq:uv}) is violated.
\end{proof}

{\bf Proof of Theorem \ref{thm:main}.} Assume the contrary:
$\Phi|_{H}$ is irreducible. Without loss we may assume that $\Phi$
is absolutely irreducible and that $H$ is a maximal subgroup of $G$.
Now we can apply the Reduction Theorem \ref{reduction} to see that
$H$ is $G$-conjugate to $P_a$ or $P_b$. However, this is a
contradition to Theorem \ref{thm:resparab}. \hfill $\Box$


\begin{thebibliography}{ABCD}

\bibitem{A}
  M. Aschbacher, {\it On the maximal subgroups of the finite classical groups},
Invent. Math. {\bf 76} $(1984)$, 469--514.

\bibitem{BK}
  J. Brundan and A. S. Kleshchev, {\it Representations of the symmetric
  group which are irreducible over subgroups}, J. reine angew. Math.
  {\bf 530} $(2001)$, 145--190.

\bibitem{Atl1}
  J. H. Conway, R. T. Curtis, S. P. Norton, R. A. Parker, and R. A. Wilson,
 {\it Atlas of Finite Groups}, Clarendon Press, Oxford, $1985$.

\bibitem{F}
 W. Feit, {\it The Representation Theory of Finite Groups},
 North-Holland Publ. Comp., Amsterdam, New York, Oxford, $1982$.

\bibitem{HH1}
  F. Himstedt and S.-c. Huang, {\it Character table of a Borel subgroup
  of the Ree groups ${^2F_4(q^2)}$}, LMS J. Comput. Math. {\bf 12} $(2009)$,
  1--53.

\bibitem{HH2}
  F.~Himstedt and S.-c.~Huang, {\it Character tables of the maximal
  parabolic subgroups of the Ree groups $\ta F_4(q^2)$}, (submitted).

\bibitem{I}
  I. M. Isaacs, {\it Character Theory of Finite Groups}, Dover
  Publications, New York, $1994$.

\bibitem{Atl2}
  C. Jansen, K. Lux, R. Parker, and R. Wilson, {\it An Atlas of Brauer
  Characters}, Clarendon Press, Oxford, $1995$.

\bibitem{KS}
  A. S. Kleshchev and J. Sheth, {\it Representations of the alternating
  group which are irreducible over subgroups}, Proc. London Math.
  Soc. (3) {\bf 84} $(2002)$, 194--212.

\bibitem{KT1}
  A. S. Kleshchev and P. H. Tiep, {\it On restrictions of modular spin
  representations of symmetric and alternating groups}, Trans. Amer. Math. Soc.
  {\bf 356} $(2003)$, 1971--1999.

\bibitem{KT2}
  A. S. Kleshchev and P. H. Tiep, {\it Representations of the general
  linear groups which are irreducible over subgroups}, Amer. J. Math.
  (to appear).

\bibitem{M}
  G. Malle, The maximal subgroups of $\ta F_4(q^2)$, {\it J. Algebra} {\bf
  139} $(1991)$, $52-69$.

\bibitem{N}
  H. N. Nguyen, {\it Irreducible restrictions of Brauer characters of the
  Chevalley groups $G_2(q)$ to its proper subgroups}, J. Algebra {\bf 320}
  $(2008)$, 1364--1390.

\bibitem{NTH}
  H. N. Nguyen and P. H. Tiep, with an Appendix by F. Himstedt,
  {\it Cross characteristic representations of $\tb D_4(q)$ are reducible
  over proper subgroups}, J. Group Theory {\bf 11} $(2008)$, 657--668.

\bibitem{Sx}
  J. Saxl, {\it The complex characters of the symmetric groups that
  remain irreducible in subgroups}, J. Algebra {\bf 111}
  $(1987)$, 210--219.

\bibitem{Sc}
  L. L. Scott, {\it Representations in characteristic $p$}, in:
  The Santa Cruz Conference on Finite Groups (Univ. California,
  Santa Cruz, Calif., 1979), pp. 319--331, Proc. Sympos. Pure Math.,
  {\bf 37}, Amer. Math. Soc., Providence, R.I., 1980.

\bibitem{Sh}
  K. Shinoda, {\it A characterization of odd order extensions of the Ree
  groups $\sp{2}F\sb{4}(q)$}, J. Fac. Sci. Univ. Tokyo Sect. I A
    Math. {\bf 22} $(1975)$, 79--102.

\bibitem{T}
  P. H. Tiep, {\it Finite groups admitting grassmannian $4$-designs},
  J. Algebra {\bf 306} $(2006)$, 227--243.

\end{thebibliography}
\end{document}